\documentclass[11pt,reqno]{amsart}
\usepackage{amsmath,amsthm,verbatim,amssymb,titletoc,enumerate}
\usepackage{bbm}
\usepackage{hyperref}
\hypersetup{colorlinks}

\newcommand{\arxiv}[1]{{\tt \href{http://arxiv.org/abs/#1}{arXiv:#1}}}

\newcommand{\floor}[1]{\left\lfloor {#1} \right\rfloor}

\newcommand{\old}[1]{}
\newcommand{\moniker}[1]{{\em (#1)}}

\newcommand{\eqindist}{\stackrel{d}{=}}
\newcommand{\red}[1]{#1}

\DeclareRobustCommand{\SkipTocEntry}[5]{}

\newtheorem{theorem}{Theorem}
\newtheorem{prop}[theorem]{Proposition}

\newtheorem{lemma}{Lemma}
\newtheorem{corollary}[theorem]{Corollary}
\newtheorem{conjecture}[theorem]{Conjecture}

\theoremstyle{remark}
\newtheorem*{remark}{Remark}

\numberwithin{counter}{section}

\theoremstyle{definition}

\def\00{\mathbf{0}}

\def\N{\mathbb{N}}
\def\Z{\mathbb{Z}}

\def\PP{\mathbb{P}}
\def\Pr{\mathrm{P}}
\def\eps{\epsilon}
\def\calB{\mathcal{B}}

\newcommand{\Sta}{{\tt S}}
\newcommand{\sleep}{{\tt s}} 
\newcommand{\fire}{{\tt F}}

\newcommand{\ARW}{{\tt ARW}}
\newcommand{\IDLA}{{\tt IDLA}}
\newcommand{\tmix}{t_{\mathrm{mix}}}
\newcommand{\Tfill}{T_{\mathrm{full}}}

\newcommand{\tfill}{t_{\mathrm{full}}}

\newcommand{\add}{{\tt A}}

\begin{document}

  \author{Lionel Levine and Feng Liang}

  \title[Exact sampling and fast mixing of Activated Random Walk]{Exact sampling and fast mixing of Activated Random Walk}

  \begin{abstract}
  	
	Activated Random Walk (ARW) is an interacting particle system on the $d$-dimensional lattice $\Z^d$. On a finite subset $V \subset \Z^d$ it defines a Markov chain on $\{0,1\}^V$. We prove that when $V$ is a Euclidean ball intersected with $\Z^d$, the mixing time of the ARW Markov chain is at most $1+o(1)$ times the volume of the ball. The proof uses an exact sampling algorithm for the stationary distribution, a coupling with internal DLA, and an upper bound on the time when internal DLA fills the entire ball. We conjecture cutoff at time $\zeta$ times the volume of the ball, where $\zeta<1$ is the limiting density of the stationary state.
	
  \end{abstract}
  
\address{Lionel Levine, Department of Mathematics, Cornell University, Ithaca, NY 14853. {\tt \url{https://pi.math.cornell.edu/~levine}}}

\address{Feng Liang, Department of Mathematics, Cornell University, Ithaca, NY 14853.}

\thanks{LL was partially supported by a Simons Fellowship, IAS Von Neumann Felowship, and NSF grant \href{https://www.nsf.gov/awardsearch/showAward?AWD_ID=1455272}{DMS-1455272}.}

\date{June 26, 2024}


\maketitle

{\centering Dedicated to the memory of Vladas Sidoravicius.\par}

  \tableofcontents
 
  \section{Introduction: Activated Random Walk}

A key feature of many complex systems is the release of stress sudden bursts. An example is the  pressure between continental plates, released in earthquakes. Bak, Tang, and Wiesenfeld called these ``self-organized critical'' (SOC) systems, and proposed a mathematial model of them, the abelian sandpile \cite{BTW,Dhar90}. But the abelian sandpile is non-universal: Even in the limit of large system size, its behavior depends delicately on the underlying graph \cite{LPS2} and on the initial condition \cite{FLW, JJ}.

One of the best candidates for a \emph{universal} model of SOC is Activated Random Walk (ARW) \cite{RS11}.
This is an interacting particle system with two species, \textbf{active} and \textbf{sleeping}.  Active particles perform random walks and fall asleep at a fixed rate $\lambda$. Sleeping particles do not move, but become active when an active particle encounters them. To make explicit the connection to SOC, sleeping particles represent stress in the system, and a single active particle can cause a burst of activity by waking up many sleeping particles.

So far, one universality result has been proved for ARW: Rolla, Sidoravicius, and Zindy \cite{RSZ} show that there is a critical mean $\zeta_c = \zeta_c(\lambda,d)$ such that for any translation-invariant and ergodic configuration $s$ of active particles in $\Z^d$ with mean $\zeta$
	\begin{equation} \label{eq:sharp} \Pr (s \text{ stabilizes}) = \begin{cases} 1, & \zeta<\zeta_c \\ 0, & \zeta>\zeta_c. \end{cases} \end{equation}
Still missing is a rigorous connection between these infinite ARW systems and their finite counterparts. For instance, we expect that the ARW stationary distribution $\mu_V$ on a finite subset $V\subset \Z^d$ has an infinite-volume limit $\mu$, and that its mean equals $\zeta_c$. We also expect the microstructure of finite ARW clusters (such as the cluster of sleeping particles formed by stabilizing $n$ chips at the origin in $\Z^d$, studied in \cite{LS1}) to converge to $\mu$ as $n \to \infty$. These conjectures are detailed in \cite{LS2}.  

Recent work has succeeded in showing that 
$\zeta_c$ is strictly between $0$ and $1$ on many transitive graphs \cite{ST18, Taggi2019, HRR, Hu, FG, JMT}, culminating in the proof by Asselah, Forien, and Gaudilli\`{e}re that $\zeta_c <1$ on $\Z^d$ for all dimensions $d$ and all sleep rates $\lambda$ \cite{AFG};
that $\zeta_c$ is continuous and strictly increasing in the sleep rate \cite{Taggi2020}, 
and tends to zero as the square root of the sleep rate in dimension $1$ \cite{ARS}; 
and that ARW at sufficiently high density takes exponential time to stabilize on a cycle \cite{cycle}.
Despite all this progress, very little is known about the behavior of Activated Random Walk at criticality. 

In this paper we examine ARW from a different perspective, by driving a finite ARW system to a stationary state. We give an exact sampling algorithm for the stationary state, and upper bound its mixing time (the time it takes to reach the stationary state from an arbitrary intitial state). \textbf{Fast mixing} is evidence for universality in that the system forgets its initial state quickly.

To see how mixing relates to universality, we can contrast ARW to the \emph{non-universal} abelian sandpile model. In contrast to \eqref{eq:sharp}, the abelian sandpile has an interval of critical means \cite{FMR}, and the problem of whether a 
sandpile on $\Z^d$ stabilizes almost surely 
is not even known to be decidable \cite{divisible}. The root cause of this non-universality is slow mixing: For example, the sandpile mixing time on both the ball $B(0,n) \cap \Z^d$ and on the torus $\mathbb{Z}_n^d$ is of order $n^d \log n$ \cite{HJL,HS}. 
This extra log factor is responsible for the non-universality of the sandpile threshold state \cite{threshold}.
Our upper bounds (Corollaries \ref{c.ball} and \ref{c.torus}), show that the extra log factor is \emph{not} present for ARW,  providing further evidence for universality.
  
\subsection{The ARW process}

Let $P$ be the transition matrix of a discrete time Markov chain on a finite state space $V \cup \{z\}$. Here $z$ is an absorbing state ($P(z,z)=1$) called the ``sink''.  Assume that every $v_0 \in V$ can access the sink, in that there exists a path $v_0, v_1, \ldots, v_k = z$ such that $P(v_{i-1},v_{i})>0$ for all $i=1,\ldots,k$. 

Given $P$ and a vector $\lambda = (\lambda_v)_{v \in V}$ with each $\lambda_v \in [0,\infty]$, we define the \textbf{ARW process} for integer $t \geq 0$ by
	\begin{equation} \label{eq:ARWproc} \sigma_{t} = \Sta [\sigma_{t-1} + \delta_{u_{t}}]. \end{equation}
In words, the state $\sigma_t$ at each discrete time step is obtained from the previous state $\sigma_{t-1}$ by adding a single active particle at $u_t$ and then stabilizing. Here
		\begin{itemize} 
	\item $\sigma_t$ takes values in the hypercube $\{0,\sleep \}^V$. The symbol $\sleep$ stands for ``sleeping'' and will be explained below.
	\item  $u_1, u_2 \ldots \in V$ is a (possibly random) sequence of vertices, called the \textbf{driving sequence};
	\item  $\delta_v(w) = 1$ if $w=v$ and $0$ otherwise;
	\item  $\Sta$ is the \textbf{stabilization operator for activated random walk} with sleep rate $\lambda$ and base chain $P$, 
which we now define.
	\end{itemize}
Following \cite{rolla2020}, consider the total ordering on $\N \cup \{\sleep \}$
	\[ 0 < \sleep < 1 < 2 < \ldots. \]
Extend addition to a commutative operation on $\N \cup \{\sleep\}$ by declaring
	\[ 0 + \sleep = \sleep \]
and
	\[ n + \sleep = n + 1 \]
for all $n \neq 0$. In particular, $s+s = s+1 = 1+s = 2$. Note that if $\sigma$ takes values in $\{0,\sleep\}$, then $\sigma + \delta_{v}$ takes values in $\{0,\sleep,1,2\}$. 

An \textbf{ARW configuration} is a map 
	\[ \sigma : V \to \N \cup \{\sleep\}. \]
If $\sigma(v) = n \geq 1$ then we say there are $n$ active particles at $v$; if $\sigma(v) = \sleep$ then we say there is one sleeping particle at $v$; and if $\sigma(v)=0$ then we say there are no particles at $v$.

A configuration taking values in $\{0,\sleep\}$ is called a \textbf{sleeping configuration}.
The stabilization operator $\Sta$ takes an arbitrary configuration $\sigma$ as input, and outputs a sleeping configuration. If $\sigma$ takes values in $\{0,\sleep\}$, then we define $\Sta [\sigma] = \sigma$. Otherwise, we obtain $\Sta [\sigma]$ by a sequence of \textbf{firings} $\fire_v$ of vertices with at least one active particle.
Firing vertex $v$ is defined in two cases, depending whether there are at least two active particles ($\sigma(v) \geq 2$) or only one ($\sigma(v)=1$). 	
	\begin{itemize}
	\item Suppose $\sigma(v) \geq 2$. To fire $v$, move one particle from $v$ to a random vertex drawn from $P(v,\cdot)$. Formally,
	\[ \fire_v [\sigma] = \sigma - \delta_v + \delta_w \text{ with probability } P(v,w) \text{ for each } w \in V \cup \{z\}. \]
	\item Suppose $\sigma(v) = 1$. To fire $v$, put the particle at $v$ to sleep with probability 
		\[ q_v := \frac{\lambda_v}{1+\lambda_v}. \] 
Otherwise move the particle from $v$ to a random neighbor drawn from $P(v,\cdot)$.  Formally,  
\[ \fire_v [\sigma] = \begin{cases} 
	 	\sigma - \delta_v + \sleep \delta_v & \text{ with probability } q_v \\
	        \sigma - \delta_v + \delta_w & \text{ with probability } (1-q_v)P(v,w)
	        \end{cases} \]
for each $w \in V  \cup \{z\}$.	        
	\end{itemize}

The domain of $\sigma$ is $V$, not $V \cup \{z\}$; so in the case $w=z$ the term $\delta_z$ is zero. This case represents a particle falling into the sink, where it is removed from the system.

We make three remarks:

\begin{enumerate}[1.]

\item Any ARW configuration $\sigma$ reaches a sleeping configuration after some (random, but almost surely finite) number of firings of active vertices.
 
This follows from our assumption in the very beginning, that $V$ is finite and every vertex can access the sink: If any particle is still active, then try to move it along a path to the sink. If it falls asleep or strays from the chosen path, then pick another active particle and try again.  There is a positive number $\eps$ such that each such trial has probability at least $\eps$ of depositing a particle in the sink. Since the number of particles is finite and cannot increase, we reach a sleeping configuration after a finite number of firings. \medskip

\item We define the \textbf{stabilization} $\Sta [\sigma]$ as the final, sleeping configuration of particles. A crucial fact is that the stabilization does not depend on the order of firings. This \textbf{abelian property}, which is proved for ARW in \cite{RS11}, also holds for a more general class of particle systems, the abelian networks \cite{BL1}. \medskip

\item The case $\lambda_v = \infty$ for all $v$ is of special interest. It is called internal DLA (\textbf{IDLA}). Note that in this case $q_v = 1$, so that each site $v$ ``absorbs'' the first particle it receives. So this process has a simple description: Each active particle moves according to the Markov chain $P$, until reaching an unoccupied site or the sink, where it remains forever.  
A fundamental quantity associated to IDLA is the time when $V$ becomes full:
	\[ \Tfill = \min \{t \,:\, \sigma_t(v) >0 \text{ for all } v \in V\}. \]

\end{enumerate}

\subsection{Plan of the paper}
Our main goal is to upper bound the mixing time of the ARW process.

We will first give a method for exact sampling from its stationary distribution (Theorem~\ref{t.exact}) and then show that the time $\Tfill$ for IDLA to fill $V$ is a \textbf{strong stationary time} for the ARW process (Theorem~\ref{t.sst}).  To upper bound the mixing time,  
it therefore suffices to upper bound $\Tfill$.  Despite the exponential size of its state space $\{0,\sleep\}^V$, we will prove that the mixing time of the ARW process is not much larger than $\#V$ (Theorem~\ref{t.mixing}). 

These three theorems are proved in Section~\ref{s.main} for a general base chain $P$. Then in Section~\ref{s.ball}, we examine the case that $P$ is simple random walk on a Euclidean ball intersected with the $d$-dimensional integer lattice $\Z^d$, with sink at the boundary of the ball. We prove that with high probability
	 \[ \Tfill \leq \# V + (\# V)^\alpha \] 
for some $\alpha<1$ (Theorem~\ref{t.fill}). 

After a brief discussion of lower bounds, we conclude with two conjectures in Section~\ref{s.conj}.

\subsection{Instructions; Resampling; Abelian Property}

The purpose of this section is to spell out the meaning of the phrase ``order of firings'' in the abelian property, so that we can give careful proofs of our main results. 

Quench the randomness of ARW into a collection of \textbf{instructions} 
$(\rho_{n,v})_{n \in \N, \, v \in V}$. The instruction $\rho_{n,v}$ dictates what will happen the $n$th time $v$ is fired: a particle at $v$ either tries to fall asleep with probability $\lambda_v / (1+\lambda_v)$, or it steps to $w$ with probability $P(v,w)/(1+\lambda_v)$.
We assume that each sequence $(\rho_{n,v})_{n \in \N}$ is independent and identically distributed (i.i.d.), and that all $\rho_{n,v}$ are independent.

Fix an ARW configuration $\sigma$, and let $(v) = (v_1, v_2, \ldots, v_m)$ be a sequence of vertices to be fired in order.  
We say that $(v)$ is a \textbf{legal execution} for $\sigma$ if $\fire_{v_{k-1}} \ldots \fire_{v_1} [\sigma] (v_k) \geq 1$ for all $k=1,\ldots,m$.  A legal execution is called \textbf{complete} for $\sigma$ if $\fire_{v_m} \cdots \fire_{v_1} [\sigma]$ is a sleeping configuration.  The \textbf{odometer} of an execution $(v)$ is the function
	\[ f(w) = \# \{1 \leq k \leq m \,:\, v_k=w \}, \]
which counts how many times each vertex $w$ is fired.  If $(v)$ is any legal execution for $\sigma$ and $(v')$ is any complete execution for $\sigma$, then the odometer of $(v)$ is less than or equal to the odometer of $(v')$ (see \cite{RS11} or \cite[Lemma 3.4]{BL1}).
 This inequality holds pointwise, both in $V$ and in the quenched instructions.  
 In particular, any two legal complete executions for $\sigma$ have the same odometer.  Therefore they use the same subset of the quenched instructions, so they result in the same sleeping configuration $\Sta [\sigma]$.  Occasionally we will write this as $\Sta_\rho [\sigma]$ to make explicit the dependence on the instructions $\rho$. 

Implicit in the definition of the ARW process is that new independent instructions are used to stabilize at each time step. But if we wish to use the previous paragraph, then the randomness for the entire chain should be expressed in terms of a \emph{single} collection of instructions $\rho = (\rho_{v,j})_{v \in V, \, j \in \N}$.   
 Our first order of business is to check that doing so does not change the distribution of $(\sigma_t)_{t \in \N}$. For this purpose we will use a lemma from \cite{LS1}.
  
For $f : V \to \N$, write $\mathcal{F}_{f}$ for the $\sigma$-field generated by the instructions $\rho_f := (\rho_{v,n})_{v \in V, n \leq f(v)}$ (the ``past''), and write $\rho^f := (\rho_{v, \, k+1+f(v)})_{v \in V, \, k\in \N}$ (the ``future'').
 
\begin{lemma} \moniker{Strong Markov Property For Quenched Instructions, \cite[Proposition 4]{LS1}}
\label{l.smp}
Let $F: V \to \N$ be a random function satisfying $\{F=f\} \in \mathcal{F}_f$ for all $f: V \to \N$.
Then $\rho^F$ has the same distribution as $\rho$, and $\rho^F$ is independent of $\rho_F$.
\end{lemma}

Now using a single collection of instructions $\rho$, let $F_t(v)$ be the number of instructions used at $v$ during the first $t$ time steps of the ARW process, let $\rho_t = \rho_{F_t}$ and let $\rho^t = \rho^{F_t}$; formally, for each $t \geq 1$ we define these inductively by 
	\begin{equation} \label{eq:ARWprecise} \sigma_t := \Sta_{\rho^{t-1}} [\sigma_{t-1} + \delta_{u_t} ] \end{equation}
and $F_t := F_{t-1} + G_t$, where $G_t$ is the odometer for stabilizing the right side. 

\begin{lemma} \label{l.resampling}
\moniker{Resampling Future Instructions}
Assume the driving sequence $u$ is independent of the instructions $\rho$.

Let $\widetilde{\rho}_1, \widetilde{\rho}_2, \ldots$ be independent families of instructions with the same distribution as $\rho$, and independent of $u$.
Then $(u_t, F_t, \sigma_t)_{t \in \N}$ has the same distribution as $(u_t, \widetilde{F}_t,\widetilde{\sigma}_t)_{t \in \N}$, where
	\[ \widetilde{\sigma}_t = \Sta_{\widetilde{\rho}_{t}} [ \widetilde{\sigma}_{t-1} + \delta_{u_t} ] \]
and $\widetilde{F}_t - \widetilde{F}_{t-1}$ is the odometer for stabilizing the right side. 
\end{lemma}

\begin{proof}
Fix functions $f_1 \leq \cdots \leq f_t$ and ARW configurations $\tau_1, \ldots, \tau_t$. 
Let $A_t$ and $\widetilde{A}_t$ be the events 
$\{ F_s = f_s, \sigma_s = \tau_s, \, s=1,\ldots,t \}$ and $\{ \widetilde{F}_s = f_s, \widetilde{\sigma}_s = \tau_s, \, s=1,\ldots,t \} $ respectively.

Write $\Pr_u$ for the law of $u$, and $\Pr$ for the law of $(\rho, \tilde{\rho}_1, \tilde{\rho}_2, \ldots)$. 
For any \emph{fixed} driving sequence $u_1, \ldots, u_t$, writing $\xi_s = \tau_{s-1} + \delta_{u_s}$ we have
	\begin{align*} \Pr (A_t) = \prod_{s=1}^t \Pr (A_s | A_{s-1}) 
			&= \prod_{s=1}^t \Pr ( \Sta_{\rho^{s-1}} [\xi_s ] = \tau_s, F_{s}=f_s | A_{s-1} )  \\
			&= \prod_{s=1}^t \Pr ( \Sta_{\rho^{s-1}} [ \xi_s ] = \tau_s , F_s=f_s) \\
			&= \prod_{s=1}^t \Pr ( \Sta_{\widetilde{\rho}_{s}} [ \xi_s ] = \tau_s , \widetilde{F}_s=f_s) \\
			&= \prod_{s=1}^t \Pr ( \Sta_{\widetilde{\rho}_{s}} [\xi_s ] = \tau_s, \widetilde{F}_s=f_s | \widetilde{A}_{s-1} ) \\
			&= \prod_{s=1}^t \Pr (\widetilde{A}_s | \widetilde{A}_{s-1}) = \Pr (\widetilde{A}_t).
		\end{align*}
In the second line we have used that the event $A_{s-1}$ depends only on the past $\rho_{s-1}$, which is independent of the future $\rho^{s-1}$ by Lemma~\ref{l.smp}. 
 In the third line we have used that $\rho^{s-1}$ has the same distribution as $\widetilde{\rho}_{s}$, again by Lemma~\ref{l.smp}.
In the fourth line we have used that $\widetilde{A}_{s-1}$ depends only on the instructions $\tilde{\rho}_r$ for $r \leq s-1$, which are independent of $\tilde{\rho}_{s}$ by hypothesis.

Now let $B = \{u_1 = v_1, \ldots, u_t=v_t\}$.
Since $u$ is assumed independent of $\rho$ and $\widetilde{\rho}$, the proof is finished by multiplying by $1_B$, taking $\Pr_u$ of both sides, and applying Fubini's theorem:
	\[ \mathbb{P} (A_t \cap B) = \Pr_u ( \Pr (A_t) 1_B ) = \Pr_u (\Pr(\widetilde{A}_t) 1_B) = \mathbb{P}(\widetilde{A_t} \cap B). \qedhere. \]
\end{proof}

As a consequence of Lemma~\ref{l.resampling}, if the driving sequence $u$ is an i.i.d.\ sequence then the ARW process \eqref{eq:ARWprecise} is a time-homogeneous Markov chain. 
For general $u$, the ARW process is not a Markov chain, but we will see that some techniques from the theory of Markov chains, such as the use of a strong stationary time to bound the mixing time, can still be applied.  The reason we are interested in general driving is for future applications when $V$ is a subset of a larger system $V'$, and the driving comes from particles entering $V$ as a result of stabilizing $V' \backslash V$.

In what follows we write $\Sta^F := \Sta_{\rho^F}$. 

\begin{lemma} \moniker{Abelian Property} \label{l.abelian}
Let $\phi_1, \phi_2$ be ARW configurations, and let $F$ be the odometer for stabilizing $\phi_1$. Then 
	\begin{equation} \label{eq:abelian} \Sta [\phi_1 + \phi_2] = \Sta^F [ \Sta[\phi_1] + \phi_2]. \end{equation}
\end{lemma}

\begin{proof}
Let $(v)$ be a legal complete execution for $\phi_1$ with instructions $\rho$, and let $(w)$ be a legal complete execution for $\Sta[\phi_1] + \phi_2$ with instructions $\rho^F$.  Then the concatenaion $(v,w)$ is a legal complete execution for $\phi_1 + \phi_2$ with instructions $\rho$.
\end{proof}

Our main use of the abelian property will be to stabilize the driving particles all at once, instead of one at a time:
	 \[ \sigma_t = \Sta [ \sigma_0 + \phi_t ] \]
where 
	\begin{equation} \label{eq:alldriving} \phi_t = \delta_{u_1} + \ldots + \delta_{u_t}. \end{equation}
This will allow us to couple the ARW and IDLA processes.

\subsection{Coupling ARW and IDLA}

Recall that ARW with infinite sleep rate is called IDLA. We write $\Sta_\infty$ for IDLA stabilization \emph{without allowing any particles to fall asleep}. In other words, to perform $\Sta_\infty$, we let each active particle perform $P$-walk until reaching an unoccupied site or the sink $z$. In particular, if $\phi$ is an all active configuration, then $\Sta_\infty [\phi]$ is all active.

Every legal IDLA execution is also legal for ARW
 (since a particle moves in IDLA only when another particle is present at the same site). 
One way to stabilize an ARW configuration is therefore to perform IDLA first, and then complete the ARW stabilization:
	\begin{equation} \label{eq:IDLAfirst} \Sta [\phi] = \Sta^{G} [ \Sta_\infty [\phi] ] \end{equation}
where $G$ is the odometer for IDLA.
In particular, with $\phi_t$ given by \eqref{eq:alldriving}, we have a coupling between the IDLA process 
	\[ \eta_{t}  = \Sta_\infty [\sigma_0 + \phi_t ] \]
and the ARW process
	\begin{equation} \label{eq:thecoupling} \sigma_{t} = \Sta [ \sigma_0 + \phi_t ] = \Sta^{G_t} [\eta_t] \end{equation}
where $G_t$ is the odometer for IDLA-stabilizing $\sigma_0 + \phi_t$. 

This coupling was used by Shellef \cite{Shellef} to prove nonfixation of certain infinite ARW systems. We will use it to bound the mixing time of the ARW process.

\section{Main Results}
\label{s.main}

Now we are ready to prove our exact sampling theorem for the ARW process $\sigma_t = \Sta [ \sigma_{t-1} + \delta_{u_{t}}]$ with driving sequence $u = (u_t)_{t \in \N}$.
We make no assumption that $u$ is i.i.d.\ or even Markovian, but we will always assume that $u$ is independent of the quenched instructions.

\subsection{Exact sampling} \label{s.exact}

 Write $\bar{\Sta} := \Sta_{\bar{\rho}}$, where $\bar{\rho}$ is an independent copy of the instructions $\rho$ used to define the ARW process.

\begin{theorem} \moniker{Exact sampling from the ARW stationary distribution} \label{t.exact}
Let $\sigma_0 = \bar{\Sta}[1_V]$. 
Then for any driving sequence $u$ and all $t \geq 1$,
	\begin{equation*} \label{eq:samedist} \sigma_t \eqindist \sigma_0. \end{equation*}
\end{theorem}

\begin{proof}
For any ARW configuration $\phi$, consider stabilizing $1_V + \phi$ in two ways: If we first move the extra particles $\phi$, then they cannot fall asleep (as every $v \in V$ contains an active particle) so they all perform $P$-walk until reaching the sink. We can then stabilize $1_V$. The second way is to stabilize $1_V$, then add the extra particles $\phi$, and stabilize again.  Using \eqref{eq:IDLAfirst} and \eqref{eq:abelian},
	\begin{equation} \label{eq:twoways} \bar{\Sta}^{G} [1_V] = \bar{\Sta} [ 1_V + \phi ] = \bar{\Sta}^{H} [ \bar{\Sta}[1_V] + \phi ].
	\end{equation}
where $G$ is the odometer for IDLA-stabilizing $1_V + \phi$ to $1_V$, and $H$ is the odometer for ARW-stabilizing $1_V$. These equalities hold pointwise in $\bar{\rho}$.

Now take $\phi = \phi_t = \delta_{u_1} + \ldots + \delta_{u_t}$. By the Strong Markov Property, in equation \eqref{eq:twoways} the left side $\eqindist \sigma_0$, and the right side $\eqindist \Sta [ \sigma_0 + \phi_t] = \sigma_t$.
\end{proof}

Theorem~\ref{t.exact} identifies a stationary distribution for the ARW process. Next we give a sufficient condition for the stationary distribution to be unique.  \red{Given $v,w \in V$ we say that $v$ can \textbf{access} $w$ if there exists $j \geq 0$ such that $P^j(v,w)>0$. We say that a fixed sequence $v_1, v_2, \ldots \in V$ is \textbf{thorough} if for all $w \in V$ there exist infinitely many $t$ such that $v_t$ can access $w$.}
Note that if the base chain $P$ is irreducible, then every $v$ can access every $w$, so \emph{every} sequence is thorough. 

\begin{lemma}
\label{l.fillup}
Let $u$ be a \red{random} driving sequence, and let $\eta_t$ be the IDLA-stabilization of $\sigma_0 + \delta_{u_1} + \ldots + \delta_{u_t}$. Then
$\PP (\eta_t = 1_V \text{ eventually}) \geq \red{\PP(u \text{ is thorough})}$.
\end{lemma}

\begin{proof} 
Let $A_t = \{v \in V \,:\, \eta_t(v)= 1 \}$. If $A_t = V$, then $A_s = V$ for all $s \geq t$. Otherwise, on the event that $u$ is thorough, it happens infinitely often that $P$-walk started at $u_{t+1}$ and stopped on exiting $A_t$ has a positive probability to exit in $V\setminus A_t$, in which case $A_{t+1}$ is strictly larger than $A_t$. Hence $\PP (A_t = V \text{ eventually}) = 1$.
\end{proof}

\red{We will be interested in random driving sequences $u$ satisfying \[ \PP(u \text{ is thorough})=1. \] (Exercising our right to omit the phrase ``almost surely,'' we will call these simply \emph{thorough}.})

Let
	\[
\mathcal{R}:=\left\{
\sigma \in \{0,\sleep\}^V \biggm| \begin{array}{l}
\text{$\sigma(v)=0$ for all $v$ such that $\lambda_v = 0$, and} \\
\text{$\sigma(w)=\sleep$ for all $w$ such that $\lambda_w=\infty$} 
\end{array}
\right\}.
\]

\begin{lemma} \moniker{Recurrent ARW Configurations} \label{l.recurrent}
If the driving sequence is thorough, then 
\begin{itemize}
\item An ARW configuration $\sigma$ is recurrent if and only if $\sigma \in \mathcal{R}$; and
\item $\mathcal{R}$ is the unique communicating class of recurrent configurations.
\end{itemize}
\end{lemma}

\begin{proof}
We first check that if $\sigma_0 \in \mathcal{R}$, then $\sigma_t \in \mathcal{R}$ for all $t$. For each vertex $v$ with $\lambda_v=0$, since $\sigma_0(v)=0$ and no particle will ever fall asleep at $v$, we have $\sigma_t(v)=0$ for all $t$.
For each vertex $v$ with $\lambda_v=\infty$, since $\sigma_0(v)=\sleep$ and the last particle left at $v$ will always fall asleep there, we have $\sigma_t(v)=\sleep$ for all $t$. 

To finish the proof, we now show that if the driving sequence is thorough, then every ARW configuration $\sigma_0$ can access every $\tau \in \mathcal{R}$.

By Lemma~\ref{l.fillup} there exists $T$ such that $\eta_T = 1_V$, so $\sigma_0 + \phi_T$ has a legal IDLA execution to $1_V$. Now starting from $1_V$, for each site $v$ such that $\tau(v) = 0$, let the particle at $v$ perform $P$-walk to the sink. Each of these walks has positive probability to reach the sink before the particle falls asleep (here we use that $\tau(v)=\sleep$ for all $v$ such that $\lambda_v=\infty$, so all such $v$ are already occupied). Then let all remaining particles fall asleep immediately. This last step succeeds with probability $\prod_{} \frac{\lambda_v}{1+\lambda_v}$, where the product is over all $v$ such that $\tau(v)=\sleep$ (here we use that $\tau(v)=0$ for all $v$ such that $\lambda_v=0$, so the product is $>0$). 
If any step fails, then repeat the whole procedure from the beginning.
\end{proof}

In the case of i.i.d.\ driving, the ARW process is a Markov chain, so uniqueness of the stationary distribution follows immediately from Lemma~\ref{l.recurrent}. 
The next lemma shows uniqueness for more general driving. We write $\Pr$ for the law of the instructions, $\Pr_u$ for the law of the driving sequence, and $\PP = \Pr_u \times \Pr$ for their joint law.

\begin{lemma} \label{c.unique}
If the driving sequence $u$ is thorough, then the ARW process has a unique stationary distribution, and the stationary distribution does not depend on $u$.
\end{lemma}

\begin{proof}
By Theorem~\ref{t.exact} the configuration $\bar{\Sta}[1_V]$ is stationary and does not depend on $u$.

To show uniqueness, let $\mu$ be any stationary distribution, and let $\sigma_0 \sim \mu$.
By stationarity of $\mu$, 
and the coupling \eqref{eq:thecoupling}, we have for all $t$ and all ARW configurations $\xi$
		\begin{align*} \mu(\xi) &= \PP (\sigma_t = \xi) 
			= \PP (\Sta^{G_t}[\eta_t]= \xi).
		\end{align*}
Now for a fixed driving sequence $u$, by the Strong Markov Property, the future instructions $\rho^{G_t}$ have the same distribution as $\bar \rho$. Since $\eta_t$ depends only on the past instructions $\rho_{G_t}$, we have (pointwise in $u$)
	\begin{align*} \Pr (\Sta^{G_t}[\eta_t]= \xi) &= \Pr (\bar{\Sta}[\eta_t] = \xi) \\
			&\geq \Pr ( \bar{\Sta}[1_V] = \xi, \eta_t = 1_V ) \\
			&= \pi(\xi) \Pr (\eta_t = 1_V)
			\end{align*}
where $\pi$ is the distribution of $\bar{\Sta}[1_V]$. Taking $\Pr_u$ of both sides,
	\[ \mu(\xi) \geq \pi(\xi) \PP(\eta_t = 1_V). \]
Since $u$ is thorough, by Lemma~\ref{l.fillup} the right side converges to $\pi(\xi)$ as $t \to \infty$.
Since both $\mu$ and $\pi$ sum to $1$ we conclude that $\mu= \pi$.
\end{proof}

In the case the stationary distribution of the ARW process is unique, we denote it by $\pi = \pi_{\lambda, P}$. 
We make a few remarks.
\medskip

\begin{enumerate}
\item Theorem~\ref{t.exact} gives a reasonably fast \textbf{sampling algorithm} for $\pi_{\lambda,P}$: The time to stabilize $\Sta_{\lambda,P} [1_V]$ is upper bounded by the time to stabilize $\Sta_{0,P} [1_V]$, which is simply the time for all particles to reach the sink $z$. Writing $T_{vz}$ for the time for a $P$-walker started at $v$ to hit $z$, the time to generate a sample from $\pi_{\lambda,P}$ is therefore at most the sum of independent hitting times 
	\[ \sum_{v \in V} T_{vz}. \] 
\medskip

\item The case when the driving sequence is constant, $u_t = v$ for all $t \in \N$, is already interesting. The ARW process $(\sigma_t)_{t \in \N}$ depends on the choice of $v$, but its stationary distribution does not. One way to see this directly is to define an operator $\add_v$ that adds one chip at $v$ and then stabilizes. 
This $\add_v$ is a stochastic matrix of size $2^{\#V}$. 
Then $\add_v \add_w = \add_w \add_v$ by Lemma~\ref{l.abelian}. The stationary distribution $\pi$ is a left eigenvector of both $\add_v$ and $\add_w$. 
\medskip

\item Despite the fast sampling algorithm, many properties of the stationary distribution $\pi_{\lambda,P}$ remain mysterious. For example, in the special case that $P$ is simple random walk on a path $\{0,1,\ldots,L\}$ with sink $z=L$, experiments indicate that $\pi_{\lambda,P}$ is \textbf{hyperuniform} in that the variance of the number of sleeping particles grows sublinearly with $L$. A number of other conjectures about $\pi_{\lambda,P}$ are discussed in \cite{LS2}.
\end{enumerate}

\subsection{Strong stationary time}

Our next goal is to show that the time $\Tfill$ for IDLA to fill $V$ is a strong stationary time for the ARW process. In words, the ARW process is exactly stationary at time $\Tfill$ and all later times. 

\begin{theorem} \moniker{Strong stationary time}
\label{t.sst}
Let $\PP= \PP_{\sigma_0,u,\lambda,P}$ be the law of the ARW process $(\sigma_t)_{t \in \N}$ with initial state $\sigma_0$, thorough driving sequence $u=(u_t)_{t \in \N}$, sleep rate vector $\lambda$, and base chain $P$ on state space $V$. 
For all ARW configurations $\sigma_0, \xi \in \{0,\sleep\}^V$, and all $t \in \N$, we have
	\begin{equation}
	\label{eq:SST}
	 \PP (\sigma_{t} = \xi \,|\, \Tfill \leq t) = \pi (\xi)
	\end{equation}
where $\pi = \pi_{\lambda,P}$ is the unique stationary distribution of the ARW process.
\end{theorem}

\begin{proof}
We will use the coupling \eqref{eq:thecoupling} between the IDLA process $\eta_t$ and the ARW process $\sigma_t = \Sta^{G_t} [\eta_t]$.
For each fixed driving sequence $u$, the event  \[ \{\Tfill \leq t\} = \{\eta_t = 1_V\} \] depends only on the past instructions $\rho_{G_t}$, which are independent of the future instructions $\rho^{G_t}$ by the Strong Markov Property. So we have (pointwise in $u$)
	\begin{align*} \Pr (\sigma_t = \xi , \Tfill \leq t) 
		&= \Pr (\Sta^{G_t} [1_V] = \xi , \Tfill \leq t) \\
		&= \Pr (\Sta^{G_t} [1_V]=\xi ) \Pr(\Tfill \leq t). 
		\end{align*}
Note that $G_t$ depends on $u$, but the future instructions $\rho^{G_t}$ can be replaced with new independent instructions $\bar{\rho}$ by the Strong Markov Property, so $\Pr (\Sta^{G_t} [1_V]=\xi ) = \Pr( \bar{\Sta} [1_V] = \xi)$ does not depend on $u$, and equals $\pi(\xi)$ by Theorem~\ref{t.exact}. So
	\[ \Pr (\sigma_t = \xi , \Tfill \leq t) = \pi(\xi) \Pr (\Tfill \leq t). \]
Recalling $\PP = P_u \times \Pr$, we obtain \eqref{eq:SST} by taking $\Pr_u$ of both sides.
\end{proof}

Bristiel and Salez have refined Theorem~\ref{t.sst} to show that the strong stationary time $\Tfill$ is separation-optimal \cite[Proposition~1]{BS}.

\subsection{Upper bounds on mixing time}
Let $u$ be a thorough driving sequence. Given a fixed (deterministic) ARW configuration $\sigma_0$, write
$\mu_t$ for the resulting distribution of the ARW process $\sigma_t$ at time $t$, and for $\eps>0$ let
	\[ \tmix ( \ARW, u,\eps ) = \min \left\{t \colon \max_{\sigma_0} ||\mu_t - \pi||_{TV} \leq \eps \right\} \]
where $|| \cdot ||_{TV}$ denotes the total variation distance between proability measures.
Let
	\[ \tfill (\IDLA,u,\eps) = \min \{t \, \colon \PP_{0, u,\infty,P}(\Tfill>t) \leq \eps \} \] 
be the first time that IDLA, started from the empty initial configuration, fills $V$ with probability $ \geq 1-\eps$.

\begin{theorem} \moniker{Upper bounds on mixing} \label{t.mixing}
\red{For any base chain $P$, any thorough driving sequence $u$, any sleep rate vector $\lambda$, and any $\eps>0$},
	\begin{equation} \label{eq:mixfill} \tmix (\ARW,u,\eps) \leq \tfill (\IDLA,u,\eps). \end{equation}
If the driving sequence $(u_t)_{t \in \N}$ is independent with the uniform distribution on $V$, then 
	\begin{equation} \label{eq:coupon} \tmix (\ARW,u,\eps) \leq \# V \log \#V + \log(1/\eps) \# V. \end{equation}
Finally, if the driving sequence $u$ is a \underline{permutation} of $V$, then the ARW process is \underline{exactly stationary} at time $\#V$, so
	\[ \tmix (\ARW,u,\eps) \leq \#V. \]
\end{theorem}

\begin{proof}
For $t \geq \tfill(\IDLA,u,\eps)$ we have by Theorem~\ref{t.sst} 
	\[ \mu_t(\xi) \geq \PP (\sigma_t = \xi, t \geq \Tfill) = \pi(\xi) \PP(t \geq \Tfill) \geq (1-\eps) \pi(\xi). \]
Summing over ARW configurations $\xi$ for which $\pi(\xi) > \mu_t(\xi)$ yields
	\[ ||\pi - \mu_t ||_{TV} = \sum_\xi (\pi(\xi) - \mu_t(\xi))_+ \leq \eps \sum_\xi \pi(\xi) \leq \eps \]
which proves \eqref{eq:mixfill}.
	
The inequality (\ref{eq:coupon}) follows from a standard coupon collector bound; see, for example, \cite[Prop.\ 2.4]{LP}, which implies that for $t>\# V \log \#V + \log(1/\eps) \# V$ we have
	\[ \PP (\phi_t \geq 1_V) \geq 1-\eps. \]
On the event $\phi_t \geq 1_V$, letting all extra particles perform $P$-walk until reaching the sink yields a legal execution from $\sigma_0 + \phi_t$ to $1_V$, so the total variation distance between the laws of $\Sta [\sigma_0 + \phi_t]$ and $\Sta [1_V]$ is at most $\eps$.

Finally, if the driving sequence is a permutation of $V$, then for $t= \# V$ we have $\phi_{\# V} = 1_V$, so $\sigma_t = \Sta [\sigma_0 + 1_V] \eqindist \Sta[1_V]$ is exactly stationary by Theorem~\ref{t.exact}.
\end{proof}

In the next section we will upper bound the right side of \eqref{eq:mixfill} when $V$ is a discrete ball in $\Z^d$.

\section{Bounds for the fill time of IDLA}
 \label{s.ball}

 For $r > 0$ let $B_r = B(0,r) \cap \Z^d = \{x \in \Z^d \,:\, |x|<r \}$ be the Euclidean ball of radius $r$ intersected with $\Z^d$, viewed as a graph with nearest-neighbor adjacencies. Here $|x | := (x_1^2+\cdots+ x_d^2)^{1/2}$ denotes the Euclidean norm.  Collapse the boundary
  	\[ \partial B_r = \{y \in \Z^d \setminus B_r \colon |y-x|=1 \text{ for some } x\in B_r \} \] 
to a sink vertex.
  	
  We consider IDLA driven by simple random walk on $B_r$, in two different scenarios: \textbf{central driving} in which all particles start at $0$, and \textbf{uniform driving} in which each particle starts at an independent random location in $B_r$. 
  
  \begin{theorem} \moniker{Upper bound for the fill time of IDLA} \label{t.fill} ~ \smallskip
  
Let $\Tfill$ be the time for IDLA with either central or uniform driving to fill $B_r$, and let $N=\# B_r$. 
  	\begin{itemize}
  	\item In dimension $d=1$, for any $\alpha>\frac{1}{2}$ there is a constant $R = R(\alpha)$ such that for all $r \geq R$
  	\begin{equation*}
	\Pr \{ \Tfill > N+N^\alpha \} \leq \exp\{- c_1 r^{\alpha-\frac{1}{2}}\}.
  	\end{equation*}
  	  	
  	\item In dimension $d\geq 2$, let $\alpha=1-\frac{1}{3d}$. Then for all sufficiently large $r$
  	\begin{equation*}
	\Pr \{ \Tfill > N+N^\alpha \} \leq \exp\{- c_2 r^{1/4}\}.
  	\end{equation*}  		
	\end{itemize}
  \end{theorem}
 
These two bounds are proved in Sections~\ref{s.dimension1} and \ref{s.higherdimensions}, respectively. The exponent~$\frac12$ is optimal for $d=1$, but $1-\frac{1}{3d}$ is not optimal for $d \geq 2$. Using methods of \cite{AG1, AG2, log,log2}, it can be improved to $1-\frac{1}{d}+\delta$, at the cost of reducing $r^{1/4}$ on the right side to $r^c$ for $c=c(d,\delta)>0$; but we do not pursue this variation. The $c_1$ and $c_2$ above are absolute constants; the proof will show that \red{$c_1=\frac{1}{266}$} and \red{$c_2 = 1$} suffice.
 
 Combining Theorem~\ref{t.fill} with the bound \eqref{eq:mixfill}, we obtain an upper bound on the mixing time of the ARW process.

\begin{corollary} \moniker{Upper bound for ARW mixing on the ball} \label{c.ball} ~ \smallskip

	Let $u$ be either the central or uniform driving sequence on the ball $B_r$, let $\lambda$ be any sleep rate vector, let P be the simple random walk on $B_r$, and let $N=\# B_r$.  Then for any $\eps>0$, we have for sufficiently large $r$,
			\begin{equation*} \tmix (\ARW,u,\eps) \leq N+N^{1-\frac{1}{3d}}. \end{equation*}
\end{corollary} 

An interesting question (see Conjecture~\ref{c.cutoff}) is whether the ARW process achieves cutoff in total variation at an earlier time $\zeta N$ for some $\zeta <1$.

By covering the torus $\mathbb{Z}_n^d$ with Euclidean balls, we obtain the following corollary, proved in Section~\ref{s.torus}. 

\begin{corollary}  \moniker{Upper bound for ARW mixing on the torus} \label{c.torus} ~\smallskip

	Let $u$ be the uniform driving sequence on the discrete torus $\mathbb{Z}_n^d \setminus \{z\}$ with sink at $z$. Let $\lambda$ be any sleep rate vector, let P be the simple random walk on $\mathbb{Z}_n^d$, and let $N=n^d$. Then for any $\eps>0$ we have for sufficiently large $n$	
	\begin{equation*} \tmix (\ARW,u,\eps) \leq N+d^{1/2}N^{1-\frac{1}{3d}}. \end{equation*}
\end{corollary}

\subsection{Concentration inequalities}

To prepare for the proof of Theorem~\ref{t.fill} we recall three concentration inequalities.

\begin{enumerate}[1.]
\item \textbf{Azuma-Hoeffding inequality.} \cite[Theorem A.10]{LP}
Let $S_t$ be a martingale with bounded differences $|S_t - S_{t-1}| \leq b_t$. Then
\begin{equation} \label{eq:Hoeffding}
	\red{\mathrm{P}(S_t-S_0 \geq s) \leq \exp \left(-\frac{s^2}{2\sum_{i=1}^{t} b_i^2} \right)}.
\end{equation}

\item \textbf{Bernstein inequality.} \cite[Theorem 2.8.4]{HDP}
Let $X_1, \ldots, X_t$ be independent mean zero random variables with $|X_i| \leq 1$. Then
	\[ \mathrm{P}(|X_1 + \ldots + X_t| \geq s) \leq 2 \exp \left(- \frac{s^2}{2 (\sum_{i=1}^t \mathrm{E}X_i^2 + \frac13 s)}\right). \]
We will apply this inequality in the case $X_i = Y_i - \mathrm{E}Y_i$ where the $Y_i$ are independent Bernoulli random variables. Writing $S = Y_1 + \ldots + Y_t$ and $\mu = \mathrm{E}S$, we obtain for $s \leq \mu$
	\begin{equation} \label{eq:bernstein} \mathrm{P}(|S - \mu| \geq s) \leq 2 \exp \left(- \frac{s^2}{3\mu}\right). \end{equation}
	
\item \textbf{Time to exit a ball.} 
Consider a simple random walk in $\mathbb{Z}^d$ starting at any point in the ball $B_r$. Let $T$ be the first exit time of the walk from $B_r$. Then for sufficiently large $t$
\begin{equation} \label{eq:SRWexit}
	\mathrm{P}(T \geq t) \leq \exp \left(-\frac{t}{3(r+1)^2} \right).
\end{equation}
This follows from the fact (proved by optional stopping for the martingale $\|W_t\|^2 - t$) that $\mathrm{E}T\leq (r+1)^2$ regardless of where the simple random walk $W_t$ starts. By Markov's inequality and the strong Markov property, $\mathrm{P}(T \geq (k+1)e(r+1)^2 \,|\, T \geq ke(r+1)^2) \leq \frac{1}{e}$ for all $k \in \mathbb{N}$. Therefore $\mathrm{P}(T \geq t) \leq (\frac1e )^{\floor{t/e(r+1)^2}}$, which implies \eqref{eq:SRWexit} for sufficiently large $t$.
\end{enumerate}

  \subsection{Upper bound in dimension 1}  
  \label{s.dimension1}

Consider IDLA with $2r+n$ particles in the interval $(-r,r)$, with particles killed if they reach an endpoint $r$ or $-r$.
By the abelian property, we may assume all particles are present at the beginning instead of being added one at a time. We stabilize IDLA in discrete time steps where at each time step, one particle moves either left by $1$ or right by $1$ with probability $1/2$ each. For definiteness, we always move the leftmost active particle (recall that a particle is active in IDLA if and only if there is at least one other particle located at the same site). We keep track of the quantity
	\[ S_t = \sum_{i=1}^{2r+n} x_{i,t} \]
where $x_{i,t}$ is the location of the $i$th particle after $t$ time steps. This $S_t$ is a martingale with \red{$|S_t - S_{t-1}| \leq 1$}; it measures the total left-right ``imbalance'' of the particles at time $t$.

\begin{lemma} \label{l.balanced}
Fix $\eps>0$ and $r^{\frac12 + \eps} \leq n \leq \red{2} r$. For IDLA with $2r+n$ particles in $(-r,r)$ with either central or uniform driving, there is a constant $R=R(\eps)$ such that for all $r \geq R$
	\begin{equation} \label{eq:balanced}
	\Pr \{ \text{No particles reach } r \} \leq \exp ( -\frac{n}{\red{33} r^{1/2}} ). 
	\end{equation}
\end{lemma}

\begin{proof}
The total number of particles is $2r+n$, and at most one particle can settle at each site in $(-r,r)$. So
on the event 
		$ \calB := \{\text{No particles reach } r\} $,
at least $n$ particles exit at $-r$. Each exiting particle contributes $-r$ to the total imbalance $S_T$, where $T$ is the time of stabilization of IDLA. If the interval $(-r,r)$ is completely full at time $T$, then the total contribution of the particles inside $(-r,r)$ to $S_T$ is zero; moreover, every unoccupied site in $(-r,r)$ results in an additional particle exiting at $-r$, which can only make $S_T$ smaller. Hence
		\begin{equation*} 
		\calB \subset \{ S_T \leq -nr \}. 
		\end{equation*}
Now for any $t \geq 0$,
		\begin{equation} \label{eq:u1} 
		\Pr (\calB) \leq \Pr \{T > t \}  +  \Pr \{ T\leq t, \mathcal{B} \} 
		\end{equation}
and since $S_t = S_T$ for all $t \geq T$, 
		\begin{align} \label{eq:u2}	
		\{ T \leq t, \mathcal{B} \} &\subset \{ T \leq t, \, S_T \leq -nr \} \nonumber \\
				&= \{T \leq t, \, S_t \leq -nr \} \nonumber \\
				&\subset \{S_t - S_0 \leq \frac{-nr}{2} \} \cup \{S_0 < -\frac{nr}{2} \}.
		\end{align}
We now make our choice of $t=nr^{5/2}$.	By Azuma-Hoeffding \eqref{eq:Hoeffding},
		\begin{equation} \red{\label{eq:u3} \Pr \{S_t-S_0 \leq \frac{-nr}{2} \} 
			\leq \exp ( - \frac{(nr/2)^2}{2t} ) 
			\leq \exp (- \frac{n}{8r^{1/2}})}. \end{equation}

With central driving, $S_0 = 0$. With uniform driving, $S_0$ is a sum of $2r+n \leq \red{4}r$ independent random variables with the uniform distribution on $(-r,r)$, so by Azuma-Hoeffding 
	\begin{equation} \red{\label{eq:driving} \Pr (S_0 < -\frac{nr}{2}) \leq \exp (- \frac{(\frac12 n r)^2}{2(4r)r^2} ) \leq \exp ( - \frac{n}{32r^{1/2}} )} \end{equation}
for sufficiently large $r$. \red{In the last inequality we have used that $n \geq r^{\frac12 + \eps}$. Combining this with \eqref{eq:u2}, \eqref{eq:u3} and \eqref{eq:driving} yields $\Pr (T \leq t, \mathcal{B}) \leq 2\exp(- \frac{n}{32r^{1/2}})$.}

Finally, to bound the first term of \eqref{eq:u1}, since \red{$n \leq 2r$},
	\[ \Pr \{T > t \} \leq \sum_{i=1}^{2r+n} \Pr \{T_i > \frac{t}{\red{4}r} \} \]
where $T_i$ is the total number of steps of taken by the $i$th particle during IDLA. By the simple random walk estimate \eqref{eq:SRWexit}, the right side is at most $\red{4}r \exp( -\frac{t}{\red{12}r(r+1)^2} ) \leq \frac12 \exp (-\frac{n}{\red{13}r^{1/2}})$ for sufficiently large $r$.
Here we again use the lower bound $n \geq r^{1/2 + \eps}$.
\end{proof}

\old{%
\begin{remark}
The only place in the proof where the driving appears is in the bound \eqref{eq:driving} on $S_0$. The right side of that bound is stronger than needed: $n^2/r$ can be replaced by $n/r^{1/2}$.
Thus, Lemma~\ref{l.balanced} holds for any initial distribution of $2r+n$ particles in $(-r,r)$ satisfying
	$ \Pr \{ S_0 < -\frac{nr}{2} \} \leq \exp (-cn/r^{1/2} ). $
In particular, this condition holds for any i.i.d. driving sequence with mean $ \geq -Cr^{1/2}$.
\end{remark}
}

Now we are ready to prove Theorem~\ref{t.fill} in dimension $1$. \red{We perform IDLA with $N+N^\alpha$ particles in the interval $B_r = (-r,r)$, where $N = 2r -1$ is the size of $B_r$.}

\begin{proof}
\red{Assume first that $\alpha \leq 1$, so that Lemma~\ref{l.balanced} applies with $n=N^\alpha$.} On the event that IDLA with central driving does not fill $B_r$, either $-r$ or $r$ receives no particles. By Lemma~\ref{l.balanced} and symmetry, this event has probability $\leq 2 \exp (-\frac{n}{\red{33} r^{1/2}} )$, which completes the proof in the case of central driving.

Now consider uniform driving. On the event that IDLA with uniform driving does not fill $B_r$, there is some $X \in B_r$ not visited by any particles, so all particles leaving $(-r,X)$ must exit to the left and all particles leaving $(X,r)$ must exit to the right. For fixed $x \in (-r,r)$ let $I$ be the larger of the two intervals $(-r,x)$ and $(x,r)$. Let $N_1$ be the number of particles starting in $I$. Then $N_1$ has the binomial$(2r+n, p)$ distribution where $p=\# I / (2r)$. Since $\# I \geq r$ we have $EN_1 \geq \# I + \frac{n}{2}$. So by Azuma-Hoeffding,
	\red{\[ \Pr \{ N_1 < \# I + \frac{n}{4} \} 
	\leq \exp ( - \frac{(n/4)^2}{2(2r+n)} ) \leq \exp (- \frac{n^2}{\red{128}r} ). \]}

On the complementary event $ \{ N_1 \geq \# I + \frac{n}{4} \}$, the probability that all $N_1$ particles exit $I$ on one side is by Lemma~\ref{l.balanced} at most
	\red{\[ \exp ( - \frac{n/4}{\red{33}(\#I / 2)^{1/2}} ) \leq \exp ( - \frac{n}{\red{132}r^{1/2}} ) \] 
since $\# I \leq 2r$. }
Taking a union bound over $x \in (-r,r)$, the probability that IDLA with $2r+n$ particles does not fill $(-r,r)$ is at most \red{$\exp (- \frac{n}{\red{133} r^{1/2}} )$} for sufficiently large $r$.

\red{For $\alpha > 1$, we split the $N+N^{\alpha}$ initial particles into batches of size $2N$, and do IDLA independently for each batch. The number of batches is $\floor{\frac{N+N^\alpha}{2N}} \geq \floor{\frac12 N^{\alpha-1}} > \frac12 r^{\alpha-1}$ for sufficiently large $r$, so by independence,
	\[\Pr \{ \Tfill > N+N^\alpha \} \leq  \exp (-\frac{r^{1/2}}{133})^{\frac{1}{2}r^{\alpha - 1}}=\exp(-\frac{1}{266}r^{\alpha - \frac12}).\]}
\end{proof}


\subsection{Upper bound in higher dimensions}  
\label{s.higherdimensions}

To prove Theorem~\ref{t.fill} in dimensions $d \geq 2$ we will use the method of Lawler, Bramson, and Griffeath~\cite{shape}. In their shape theorem for IDLA, driving is from the origin and there is no sink. Here we adapt their method to uniform driving with sink.

Fix $d \geq 2$ and let $G_r(y,z)$ be the expected number of visits to $z$ by simple random walk started at $y$ before exiting the ball \red{$B_r = \{x \in \mathbb{Z}^d \,:\, |x|<r\}$}. We recall 
(
\cite[\textsection 1.5, 1.6]{Lawler-intersections} 
that $G_r$ is symmetric in $y$ and $z$, and for all $z \in B_r$
	\begin{equation} \label{eq:escape} G_r(z,z) \leq c_1 \log r \end{equation}
for a constant $c_1$ depending only on $d$. In the proof of Theorem~\ref{t.fill}, we will use the following lower bound. 

\begin{lemma} \label{l.exittime}
There is a constant $c_2 > 0$ depending only on $d$, such that for all $z \in B_r$
	\[ \sum_{y \in B_r} \mathrm{P}_y(\tau_z < \tau_r) \geq c_2 \frac{r}{\log r}. \]
\end{lemma}

\begin{proof}
Recall that $G_r(z,z) \mathrm{P}_y(\tau_z < \tau_r) = G_r(y,z) = G_r(z,y)$; the first equality follows from the strong Markov property by noting that if the walk visits $z$ before exiting $B_r$, then the number of visits to $z$ before exiting $B_r$ has a geometric distribution with mean $G_r(z,z)$. Now by \eqref{eq:escape},
	\begin{align*}
	c_1 \log r \sum_{y \in B_r} \mathrm{P}_y(\tau_z < \tau_r)  
		&\geq \sum_{y \in B_r} G_r(z,y).
	\end{align*}
The right side equals the expected time $\mathrm{E}_z \tau_r$ for simple random walk started at $z$ to exit $B_r$. As a function of $z$, this expected time has discrete Laplacian $-1$ and vanishes outside $B_r$, so 
	\[ \mathrm{E}_z \tau_r \geq r^2 - |z|^2 \geq r(r-|z|). \]
This completes the proof for all $z \in B_{r-1}$. For $z \in B_r - B_{r-1}$ note that simple random walk started at $z$ has a constant probability of hitting $B_{r-1}$ before exiting $B_r$, so $\mathrm{E}_z \tau_r$ is at least a constant times $r$.
\end{proof}

\red{
We will use the following Green function inequality to address the case of central driving. It differs in two respects from \cite[Lemma 3]{shape}, in which the $N^{1-\frac1d}$ term is absent but $z$ is restricted to the smaller ball $B_{(1-\eps)r}$.
}

\begin{lemma}
\label{l.divsand}
\red{
Fix $r$ and let $N = \# B_r$. There is a constant $C_1$ depending only on $d$, such that
	\begin{equation} \label{eq:divsand} (N+ C_1 N^{1-\frac1d}) G_r(0,z) \geq \sum_{y \in B_r} G_r(y,z) \quad \text{ for all } z \in B_{r}. \end{equation}
}
\end{lemma}

\begin{proof} 
\red{
By \cite[Theorem 3.3]{sharpcirc}, there are constants $C_0$ and $C_2$ depending only on $d$, and a nonnegative function $u$ (the ``divisible sandpile odometer function'') on $\Z^d$ such that 
	\begin{itemize}
	\item $u$ has discrete Laplacian $1 - N_0 \delta_0$ in $B_r$, where $N_0 = N+C_0 N^{1-\frac1d}$.
	\item $u \leq C_2$ on $\partial B_r$. (This is the last displayed equation in the proof of \cite[Theorem 3.3]{sharpcirc}.)
	\end{itemize}
The function
	\[f(z) := N_0 G_r(0,z) - \sum_{y \in B_r} G_r(y,z). \] 
also has discrete Laplacian $1-N_0\delta_0$ in $B_r$, and $f$ vanishes on $\partial B_r$. By the maximum principle, $f-u \geq -C_2$ in $B_r$. 
}\red{
Now we'll use the fact that $G_r(0,z) \geq C_3 N^{-(1-\frac1d)}$ in $B_r$, for a constant $C_3$ depending only on $d$. (This follows from \cite[Prop.\ 1.5.9, 1.6.7]{Lawler-intersections} when $z \in B_{r-c}$ for a suitable constant $c$; to extend to $z \in B_r \setminus B_{r-c}$, choose $z' \in B_{r-c}$ with $|z-z'|<c$ and note that simple random walk started from $z'$ has a constant probability to visit $z$ before exiting $B_r$). Letting $C_1 = C_0 + \frac{C_2}{C_3}$ we have
	\begin{align*} (N + C_1 N^{1-\frac1d}) G_r(0,z) - \sum_{y \in B_r} G_r(y,z) &\geq f(z) + \frac{C_2}{C_3} N^{1-\frac1d} G_r(0,z) 	\\
							&\geq u(z) - C_2 + C_2.
	\end{align*}
Since $u \geq 0$, this completes the proof.
}
\end{proof}

\begin{proof}[Proof of Theorem~\ref{t.fill} in dimensions $d \geq 2$]
	We consider first the case of uniform driving. Perform IDLA starting with \red{$\floor{N+N^\alpha}$} particles at independent uniform locations in the ball $B_r$, where $N=\# B_r$. Denoting by $A_r$ the resulting random subset of $B_r$ where particles stabilize, we must show that $\Pr (A_r \neq B_r) \leq \exp (- \frac14 r^{1/4})$ for sufficiently large $r$.
	
	We modify the proof of the inner bound in \cite{shape} to account for killing at $\partial B_r$ and uniform driving. For $ z \in B_{r}$, denote by $E_z = \{z \not \in A_r\}$ the event that no particle visits $z$ during IDLA. By a union bound over $z$, it suffices to show that for sufficiently large $r$
	\begin{equation}
		\red{\mathrm{P}(E_z)< \exp (- 2 r^{1/4})} \quad \text{ for all } z\in B_{r}.
	\end{equation}
	
	Fix an arbitrary ordering of the particles, and define 
	\begin{equation}
	\begin{split}
	&\tau_z^i=\text{time of first visit to $z$ by the $i$th particle in simple random walk};\\
	&\tau_r^i=\text{time of first exit of } B_r \text{ by the $i$th particle in simple random walk};\\
	&\sigma_i=\text{stopping time of $i$th particle in the IDLA stabilization process\red{;} }\\
	&\red{M=\sum_{i=1}^{\floor{N+N^\alpha}} \mathbbm{1}_{\{\tau_z^i<\tau_r^i\}};}\\
	&\red{L=\sum_{i=1}^{\floor{N+N^{\alpha}}} \mathbbm{1}_{\{\sigma_i<\tau_z^i<\tau_r^i\}};}\\
	&\tilde{L}=\sum_{y \in B_r} \mathbbm{1}_{\{\tau_z^y<\tau_r^y\}}. 
	\end{split}
	\end{equation}
Here $\tau_z^y$ is first hitting time of $z$ for a simple random walk started at $y$; and $\tau_r^y$ is the first exit time of $B_r$ for a simple random walk started at $y$.
	
	Now we have for any $a\in \mathbb{R}$,
	\begin{align} \label{eq:LBG}
		\mathrm{P}(E_z)=\mathrm{P}(M-L=0) 
			&\leq \mathrm{P}(M \leq a)+\mathrm{P}(L \geq a) \nonumber \\
		&\leq \mathrm{P}(M \leq a)+\mathrm{P}(\tilde{L} \geq a).
	\end{align}
	The last inequality follows from the observation that after IDLA stabilization, each vertex can be occupied by at most one particle, so $\tilde{L} \geq L$.

Next we will show $\mathrm{E}M$ is substantially larger than $\mathrm{E}\tilde{L}$.	
Since \red{$\floor{N+N^{\alpha}}$} particles are dropped uniformly in $B_r$, and $N=\# B_r$, we have:
	\begin{equation} \label{eq:EM}
	\mathrm{E}M=\frac{\red{\floor{N+N^{\alpha}}}}{N}\sum_{y \in B_r} \Pr (\tau_z^y < \tau_r^y) = (1+\red{\frac{\floor{N^\alpha}}{N}}) \mu\red{,}
	\end{equation}
where $\mu = \mathrm{E}\tilde{L}$. By Lemma~\ref{l.exittime} we have $\mu \geq c_2 \frac{r}{\log r}$ for all $z \in B_r$. Now we make our choice of exponent: $\alpha = 1-\frac{1}{3d}$, so that for sufficiently large $r$
	\begin{equation} \label{eq:alphachoice}  \red{\frac{\floor{N^\alpha}}{N}} \geq c_3 r^{d(\alpha-1)} = c_3 r^{-\frac13} \geq c_3 \mu^{-\frac13-\eps} \end{equation}
for any $\eps>0$ and $r \geq R(\eps)$, where the constant $c_3>0$ depends only on $d$. This implies
	\[ \mathrm{E}M - \mathrm{E}\tilde{L} \geq \red{c_3} \mu^{\frac23 - \eps}. \]
Taking $a = (\mathrm{E}\tilde{L} + \mathrm{E}M)/2$ we have
		\begin{align*}
		\label{eq:Ltail} 
		\Pr (\tilde{L} \geq a) &\red{ \leq } \Pr (\tilde{L} - \mu \geq \frac12 \red{c_3} \mu^{\frac23 - \eps} )
		\end{align*}
By Bernstein's inequality \eqref{eq:bernstein}, the right side is at most $2 \exp(-\frac{c_3^2}{12} \mu^{\frac13 - 2\eps} )$.

Likewise, since $\mu \leq \mathrm{E}M \leq 2\mu$, we have by Bernstein's inequality
		\begin{align*}
			\Pr (M \leq a) \leq \Pr (M - \mathrm{E}M \leq -\frac14 \red{c_3} (\mathrm{E}M)^{\red{\frac23 -\eps}} ) \leq 2 \exp(\red{-\frac{c_3^2}{48} \mu^{\frac13 - 2\eps}} ).
		\end{align*}	
We conclude from \eqref{eq:LBG} that \red{for sufficiently small $\epsilon$ and sufficiently large $r$}
	\begin{equation} \label{eq:holeprob}
		\red{\Pr (E_z) \leq \exp (- 2 r^{1/4} ) }
	\end{equation}
which completes the proof in the case of uniform driving.

Now we adapt the proof to handle the case of central driving. Note driving enters the proof only in equation \eqref{eq:EM}. In the case of central driving, we have instead
	\[ \mathrm{E}M = \red{\floor{N+N^\alpha}} \Pr (\red{\tau_z^o \leq \tau_r^o}) \]
\red{so that
	\begin{equation} \label{eq:meandiff} \mathrm{E}M - \mu = \frac{1}{G_r(z,z)} \left( \floor{N+N^\alpha} G_r(o,z)- \sum_{y \in B_r} G_r(y,z) \right). \end{equation}
To complete the proof in this case, we use Lemma~\ref{l.divsand} 
to lower bound \eqref{eq:meandiff} for sufficiently large $r$ as follows:
	\[  \mathrm{E}M - \mu \geq \frac{\frac12 \floor{N^\alpha}}{N+N^\alpha} \mathrm{E}M. \]
Now applying \eqref{eq:alphachoice} to lower bound the right side, we obtain for sufficiently large $r$
	\[ \mathrm{E}M - \mu \geq (\frac14 c_3 \mu^{-\frac13 - \epsilon}) \mathrm{E}M \geq \frac14 c_3 (\mathrm{E}M)^{\frac23 - \epsilon}. \]
The proof now proceeds as before: Again choosing $a = (\mathrm{E}\tilde{L} + \mathrm{E}M)/2$, we have
	\[ \Pr ( M \leq a) \leq \Pr (M - \mathrm{E}M \leq -\frac18c_3 (EM)^{\frac23 - \eps}) 
		\]
	\[ \Pr ( \tilde L \geq a )\leq \Pr (\tilde L - \mathrm{E}\tilde{L} \geq \frac18c_3 (EM)^{\frac23 - \eps}). 
		\]
By Bernstein's inequality, since $\mathrm{E}M \geq \mathrm{E}\tilde{L}  = \mu$, both of these are upper bounded by $ 2\exp ( - \frac{c_3^2}{192} \mu^{\frac13 - 2\eps})$, which yields \eqref{eq:holeprob} as before.
}
\end{proof}

\subsection{Upper bound for the torus} \label{s.torus}

\begin{proof}[Proof of Corollary~\ref{c.torus}]
The case $d=1$ is immediate from Theorem~\ref{t.fill}.

For $d \geq 2$, we cover the torus $\mathbb{Z}_n^d \setminus \{z\}$ by Euclidean balls of radius of $n/4$, while leaving the sink $z$ uncovered. A simple but inefficient way to do this, which suffices for our purpose, is to take \emph{all} balls $B(x,n/4)$ for $x \in \mathbb{Z}_n^d$ which do not contain $z$.

Now let $B$ be one of the covering balls. By the abelian property, if IDLA with sink at $\partial B$ fills $B$, then IDLA with sink at $z$ also fills $B$.

One can check that for all $d \geq 2$ we have
$( \#B / n^d)^{1/3d} \geq 1.05 d^{-1/2}$ for sufficiently large $n$.
When $d=2$, this bound follows from the fact that $\# B \geq 3 (n/4)^2$ for sufficiently large $n$ (since $\pi>3$).
For general $d$, the bound follows from the formula for the volume of the $d$-dimensional ball, along with the estimates $k^k / e^{k-1} \leq k! \leq (k+1)^{k+1} / e^k$.

Let $n$ be large enough so that the bounds in Theorem~\ref{t.fill} hold for $r=n/4$. 
Let $N=n^d$ and $\alpha = 1-\frac{1}{3d}$. After dropping $t= \red{\floor{N+d^{1/2}N^{\alpha}}}$ particles
uniformly at random in $\mathbb{Z}_n^d \setminus \{z\}$, the number of particles starting in $B$ is a sum of $t$ independent Bernoulli random variables of mean $\#B / N$.  This sum has mean $\geq (\#B) + 1.05 (\# B)^\alpha$, so by Bernstein's inequality \eqref{eq:bernstein},
the probability that $B$ starts with less than $(\# B)+ (\# B)^\alpha$ particles is at most
$2\exp(-c(\# B)^{2\alpha-1})$, \red{where $c = \frac{0.05 ^ 2}{3}$}.  For sufficiently large $n$, on the event that $B$ starts with at least $(\# B)+ (\# B)^\alpha$ particles, the probability that IDLA does not fill $B$ is at most $\exp(-c_2(n/4)^{1/4})$, by Theorem~\ref{t.fill}.

By a union bound over the covering balls, the probability that IDLA does not fill $\mathbb{Z}_n^d \setminus \{z\}$ is at most $n^d [2\exp(-c(\# B)^{2\alpha-1}) + \exp(-c_2(n/4)^{1/4})]$. Taking $n$ large enough so that this probability is $<\eps$, we obtain from \eqref{eq:mixfill}
	\[\tmix(\ARW ,u,\eps) \leq \tfill(\IDLA, u, \eps) \leq t. \qedhere \] 
\end{proof}

\subsection{Lower bounds}
In this section we state some lower bounds for the fill time of IDLA. The proofs are straightforward, so we indicate only the main idea.

The first lower bound shows that the exponent $\alpha$ in Theorem~\ref{t.fill} cannot be improved to less than $\frac12$ in dimension $d=1$ or $1-\frac{1}{d}$ in dimensions $d \geq 2$. 

\begin{prop}
For $d \geq 1$, let $\Tfill$ be the time for IDLA to fill the ball $B_r \subset \Z^d$, with sink at $\Z^d \setminus B_r$. Let $N = \# B_r$ and let 	 $\beta=\max\{\frac{1}{2},1-\frac{1}{d}\}$.
The following holds for any driving sequence $u$ satisfying $u_t\in B_{r-2}$ for all~$t$ \red{for $d \geq 2$}, and also for the uniform driving sequence on $B_r$ \red{for $d \geq 1$}: For all $b>0$ there exists $c>0$ such that for all sufficiently large $r$
	\[
	 \mathrm{P}(\Tfill > N+bN^{\beta}) > c. 
	\]
\end{prop}

\red{The proof idea for $d = 1$ is based on the observation that with some  positive probability (depending only on $b$)  the left half interval $(-r,0)$ starts with more than $r + 4b N^{1/2}$ particles.  Conditioned on this event, with at least $\frac12$ probability, at least $2b N^{1/2}$ particles exit $(-r,0)$ at $-r$, in which case there are not enough particles remaining to fill up $(-r,r)$. }

The \red{proof idea for $d \geq$ 2} is to split the IDLA stabilization into two stages: In stage one, stabilize all particles starting inside $B_{r-2}$, stopping them when they hit $\partial B_{r-2}$; and in stage two, finish the stabilization procedure. Let $M$ be the number of particles resting at $\partial B_{r-2}$ at the end of stage one. If $M$ is large ($\geq CN^{\beta}$), then with nonvanishing probability, at least $2bN^{\beta}$ particles will exit $B_r$ in stage two. If $M$ is small ($< CN^{\beta}$), then with nonvanishing probability, $B_r \setminus B_{r-2}$ will not fill up in stage two. \\

Next we observe that for general graphs, $\Tfill$ is not always upper bounded by $(1+o(1))\#V$. The wired tree provides an example: Let $V$ be the graph obtained from the complete binary tree of depth $n+1$, by collapsing all $2^{n+1}$ leaves to a single sink vertex, $z$. We show that on this graph, the fill time of IDLA has order $\#V \log \# V$.
Denote by $\mathcal{B}$ the set of $2^n$ neighbors of $z$, and let
	\[ \Tfill' = \min \{t \,:\, \mathcal{B} \subset A_t \} \]
be the first time IDLA contains $\mathcal{B}$. Note $\Tfill \geq \Tfill'$.  

\begin{prop} 
For IDLA driven by simple random walk on the wired tree $V$, with either central or uniform driving, we have that for any $c<1/4$
	\[ \mathrm{P}(\Tfill' > c \, \#V \log \#V) \to 1 \] 
as $n \to \infty$.
\end{prop}

The idea of the proof is to lower bound the time to collect $2^n$ coupons corresponding to the vertices of $\mathcal{B}$. Each time a vertex of $\mathcal{B}$ joins the IDLA cluster, we collect a coupon. The coupons are not independent, but there is a uniform upper bound on the probability of collecting a new coupon.  This uniform upper bound is derived by continuing the path of a particle that joins the cluster until it hits $z$, so that the probability of collecting a new coupon at time $t+1$ is at most the expected number of hits of $\mathcal{B}\setminus A_t$ by simple random walk before hitting $z$.  For either central or uniform driving, this expected number of hits \red{equals} $2(1-\frac{k}{2^n})$, where $k = (\# A_t \cap \mathcal{B})$ is the number of boundary vertices in the current cluster.  Since this upper bound depends on the cluster only via $k$, we can lower bound $\Tfill'$ by a sum of $2^n$ independent geometric random variables as in the standard coupon collector.

\section{Conjectures} 
\label{s.conj}

We conclude by stating two conjectures. 
\red{While this paper was under review, Conjecture~\ref{transitive graph} has been proved by Bristiel and Salez \cite[Proposition 1 and Corollary 2]{BS}. Regarding Conjecture~\ref{c.cutoff}, which remains open, the related inequality $\zeta_c < 1$ has been proved in all dimensions and for all sleep rates \cite{AFG}.}

\begin{conjecture}\label{transitive graph} \moniker{Time for IDLA to fill a transitive graph} ~ \smallskip

	Let $V$ be a transitive graph with one vertex designated as sink. Then 
		\[ \frac{\Tfill}{\# V} \to 1 \qquad \text{in probability as } \#V \to \infty. \]
\end{conjecture}

\begin{conjecture}\label{c.cutoff} \moniker{Cutoff for ARW at the stationary density}  ~ \smallskip

	Let $u$ be the uniform driving sequence on $B_r = B(0,r) \cap \Z^d$. Let $0 < \lambda < \infty$ be any constant sleep rate, and let $P$ be the simple random walk on $B_r$ with sink at $\Z^d \setminus B_r$. There exists a constant $\zeta = \zeta(\lambda,d)<1$ such that 
	
	\begin{enumerate}
		\item For any $\eps>0$,
			\[ \frac{\tmix(\ARW,u,\eps)}{\# B_r} \to \zeta \qquad \text{as } r \to \infty. \]
		
		\item Writing $|\Sta [1_{B_r}]|$ for the number of particles in the ARW stationary state on $B_r$, we have
					\[ \frac{|\Sta [1_{B_r}]|}{\# B_r} \to \zeta \qquad \text{in probability as } r \to \infty. \]

		\item $\zeta = \zeta_c$, the critical density for ARW stabilization in $\Z^d$.
	\end{enumerate}
\end{conjecture}

\section*{Acknowledgements}

The first author thanks Christopher Hoffman, Leonardo Rolla, and Vladas Sidoravicius for inspiring conversations.
The second author thanks his parents for their long time supporting him. Both authors thank the anonymous referee for detailed comments that improved the paper.

\end{document}